\documentclass[graybox]{svmult}

\usepackage{mathptmx}       % selects Times Roman as basic font
\usepackage{helvet}         % selects Helvetica as sans-serif font
\usepackage{courier}        % selects Courier as typewriter font
\usepackage{type1cm}        % activate if the above 3 fonts are
\usepackage{makeidx}         % allows index generation
\usepackage{graphicx}        % standard LaTeX graphics tool
\usepackage{multicol}        % used for the two-column index
\usepackage[bottom]{footmisc}% places footnotes at page bottom

\usepackage{amsfonts}
\usepackage{amsmath}
\usepackage{amssymb}
\usepackage{nicefrac}
\usepackage{subfig}
\usepackage{epstopdf}
\usepackage{bbm}

\newcommand{\invLa}{\nicefrac{1}{\lambda}}
\newcommand{\ind}{\mathbbmss{1}}
\newcommand{\FF}{\mathcal{F}}
\newcommand{\EW}{\mathbb{E}}

\newcommand{\hatt}{\hat{t}}

\makeindex             % used for the subject index
\begin{document}

\title*{Optimal inflow control penalizing undersupply in transport systems with uncertain demands}
\titlerunning{Optimal inflow control penalizing undersupply}
\author{Simone G\"{o}ttlich, Ralf Korn, Kerstin Lux}
\institute{Simone G\"{o}ttlich \at University of Mannheim, Department of Mathematics, 68131 Mannheim, Germany, \email{goettlich@uni-mannheim.de}
\and Ralf Korn \at TU Kaiserslautern, Department of Mathematics, P.O. Box 3049, 67653 Kaiserslautern, Germany \email{korn@mathematik.uni-kl.de} \and Kerstin Lux \at University of Mannheim, Department of Mathematics, 68131 Mannheim, Germany, \email{klux@mail.uni-mannheim.de}}
\maketitle

\abstract{	We are concerned with optimal control strategies subject to uncertain demands. An Ornstein-Uhlenbeck process describes the uncertain demand. The transport within the supply system is modeled by the linear advection equation. 
We consider different approaches to control the produced amount at a given time to meet the stochastic demand in an optimal way. In particular, we introduce an undersupply penalty and analyze its effect on the optimal output in a numerical simulation study.
}

\section{Introduction}
\label{sec:intro}

In many real-world situations, taking uncertainty into account becomes more and more important. In the context of supply chain management, a need for appropriate control strategies under uncertainty naturally arises when it comes to production planning. The size and timing of product orders is often not known in advance.
However, for a delivery on time, the production process needs to be started in advance. In this work, we tackle the challenging question of when to feed how many goods into a supply system to meet the stochastic demand.
We use the framework of ~\cite{Lux.2018}, where a corresponding stochastic optimal control problem is set up in the context of electricity injection and extend it by introducing a penalty term into the cost function. This term penalizes a production not leading to demand satisfaction, i.e. an undersupply.

The main contribution of this work is to provide insight into the effect of an undersupply penalty on the optimal production plan.
In a numerical simulation study, we highlight the effect for different penalty parameters.

\section{Stochastic optimal control model for transport systems}
\label{sec:model}
An analysis of optimal control strategies for a supply problem in a deterministic demand setting can be found in ~\cite{Teuber.2018, Schillen.2016}. 
Here, we focus on the stochastic nature of the demand and start from the stochastic optimal control framework originally set up in ~\cite{Lux.2018}.
We consider a supply system consisting of only one production line. 
Goods are fed into the system at $x=0$, and leave the system at $x=1$.
Within a finite time interval $[0,T]$, the aim is to optimally match the externally given customers' demand $Y_t$ located at $x=1$ by determining the inflow control $u(t)\in L^2$ of goods at $x=0$.
Thereby, the transport of goods $z=z(x,t)$ along the production line is governed by the linear advection equation with constant transport velocity $\lambda>0$ and the following initial and inflow conditions:
\begin{align}
	z_t + \lambda z_x &= 0, \quad x \in (0,1),\ t \in [0,T] \notag \\
	z(x,0) &= 0, \quad z(0,t) = u(t).\label{eq:transportDyn}
	%z(x=1,t) &= y(t). \label{eq:transportDyn}
\end{align}
We denote by $y(t)=z(1,t)$ the output of the system. It is intended to match the externally given demand $Y_t$.

The uncertainty about the height and timing of the orders entails the stochasticity of $Y_t$.
As in ~\cite{Lux.2018}, we assume that the demand process fluctuates around a given time-dependent mean demand level $\mu(t)$. The latter can be seen as a forecast that is based on historical demand data. In this demand setting, one possible model choice is the Ornstein-Uhlenbeck process (OUP). Let $W_t$ be a one-dimensional Brownian motion, $\sigma>0,\ \kappa>0$ be constant parameters, and denote the initial demand by $y_0$. Then, the OUP is the unique strong solution of the stochastic differential equation (SDE)
\begin{align}
dY_t=\kappa\left(\mu(t)-Y_t\right)dt + \sigma dW_t,\ \quad Y_{0}=y_0. \label{eq:OUP}
\end{align}
The OUP possesses a mean-reverting property, i.e., whenever the process is away from its mean demand level it is attracted back to it. The parameter $\kappa$ describes how strong this attraction is and $\sigma$ determines how large the fluctuations are.

%From its explicit form, one can directly derive the distribution of $Y_t$ and obtains the following normal distribution:
In this work, we make use of the known distribution of $Y_t$, which is given by the following normal distribution:
\begin{align}
Y_t \sim N\left( {y_0e^{ - \kappa t}  + \kappa \int\limits_0^t {e^{ - \kappa \left( {t - s} \right)} \mu \left( s \right)ds} ,\,\,\sigma ^2 \int\limits_0^t {e^{ - 2\kappa \left( {t - s} \right)}ds} \,} \right). \label{eq:densityOUP}
\end{align}
We refer the reader to ~\cite{Lux.2018} for more details on the demand process and the possibility to include jumps in the demand.

The problem of interest is the arising constrained stochastic optimal control (SOC) problem
\begin{align}
\min_{u(t), t \in [0,T-\nicefrac{1}{\lambda}], u \in L^2} &  \int_{\nicefrac{1}{\lambda}}^{T} OF(Y_s,t_0,y_{t_0},y(s)) ds \ \text{subject to} \ \eqref{eq:transportDyn} \ \text{and} \ \eqref{eq:OUP} , \label{eq:SOC}
\end{align}
Thereby, $\nicefrac{1}{\lambda}$ is the time that one good needs to pass the production line, and $OF(Y_s,t_0,y_{t_0},y(s))$ denotes the loss function.

In ~\cite{Lux.2018}, a possible choice of an objective function as a tracking-type function $OF_{\text{track}}(Y_s,t_0,y_{t_0},y(s)) = \EW\left[(Y_s - y(s))^2 |Y_{t_0} = y_{t_0} \right]$ has been introduced. The loss is measured in terms of the quadratic deviation between output at the end of the line and the actual demand. In this work, we focus on an extended loss quantification including an undersupply penalty.
This is of interest for companies where a supply guarantee is of crucial importance and short-term external purchase is very costly. For them, it might be more harmful to generate an output that does not lead to demand satisfaction compared to an overproduction. Therefore, we introduce a new term into the objective function that penalizes undersupply. Thereby, $\alpha$ regulates the intensity of penalization.
\begin{align}
	OF_{\text{pen}}(Y_s,t_0,y_{t_0},y(s)) = &\EW\left[(Y_s - y(s))^2 |Y_{t_0} = y_{t_0} \right] \notag\\
	&+ \alpha \EW\left[(Y_s - y(s))^2 | Y_s > y(s) \wedge Y_{t_0} = y_{t_0}\right]. \label{eq:OFpen}
\end{align}
According to ~\cite[Def.\ 8.9]{Klenke.2008}, the second conditional expectation in \eqref{eq:OFpen} reads as
\begin{align*}
	&\EW\left[(Y_s - y(s))^2 | Y_s > y(s) \wedge Y_{t_0} = y_{t_0}\right] \\
	=& \left\{\begin{array}{ll}
	 \frac{\EW\left[(Y_s - y(s))^2 \ind_{\{Y_s > y(s)\}} | Y_{t_0} = y_{t_0}\right]}{P\left(Y_s>y(s)\right)} & \text{if} \ P(Y_s > y(s))>0 \\
	0 & \text{else}.
	\end{array}\right.
\end{align*}
Thus, both conditional expectations in \eqref{eq:OFpen} can be expressed in terms of the known demand density $\rho_{Y_t|Y_{t_0}=y_{t_0}}$ at time $t$ given by \eqref{eq:densityOUP}. Hence, for the evaluation of the objective functions $OF_{\text{track}}$ and $OF_{\text{pen}}$, this information on the demand density is sufficient. As the objective function is the only part of the SOC problem where the stochastic demand dynamics \eqref{eq:OUP} come into play, we can replace the SDE constraint \eqref{eq:OUP} in \eqref{eq:SOC} by the condition that $Y_t$ has demand density \eqref{eq:densityOUP}, which is used to calculate the expectations in the objective function \eqref{eq:OFpen}. We are left with
\begin{align}
\min_{u(t), t \in [0,T-\nicefrac{1}{\lambda}], u \in L^2} &  \int_{\nicefrac{1}{\lambda}}^{T} OF_{\text{pen}}(Y_s,t_0,y_{t_0},y(s)) ds \ \text{subject to} \ \eqref{eq:transportDyn} \ \text{and} \ \eqref{eq:densityOUP}. \label{eq:SOCmod}
\end{align}
We are now able to apply deterministic optimization algorithms to the SOC problem \eqref{eq:SOCmod}.

However, we still need to make assumptions on the demand information that is used to determine the optimal inflow $u(t)$. Those assumptions result in different control methods due to the measurability assumptions on the inflow control $u(t)$. We focus on two of the three presented control methods (CM) in ~\cite{Lux.2018} corresponding to two information scenarios that are shortly summarized here for the sake of completeness:
\begin{itemize}
	\item \textbf{CM1}: The only available demand information is the initial demand $y_0$ and the demand dynamics \eqref{eq:OUP}. No updates on the actual evolution of the demand can be used to determine the inflow control over the optimization horizon $[0,T]$. Thus, we assume that $u(t)$ is $\FF_t$-measurable, where $\FF_t = \sigma\left(Y_s; 0\leq s \leq t\right)$. 
	\item \textbf{CM2}: We prespecify update times $0=\hatt_0<\hatt_1<\cdots<\hatt_n\le T-\nicefrac{1}{\lambda}$, where $\hatt_i=i \cdot \Delta t_{\text{up}}$, $i \in \{0,1,\cdots,\nicefrac{T-\invLa}{\Delta t_{\text{up}}}\}$, and update frequency $\Delta t_{\text{up}} \in [0,T-\invLa]$. At those points in time, the initial demand and the demand dynamics \eqref{eq:OUP} are supplemented by the actually realized demand. The forecast is updated accordingly and the optimal inflow control is calculated based on the updated demand forecast. Hence, we assume $u(t)$ is $\FF_{\hatt_i}$-measurable for $t \in [\hatt_i,\hatt_{i+1}]$.
\end{itemize}

CM1 is directly applicable to \eqref{eq:SOCmod}. For CM2, we divide the optimization period $[0,T]$ into smaller subperiods $[\hatt_i,\hatt_{i+1}]$ according to the prespecified update times $\hatt_i$ and solve our SOC problem thereon.
\begin{align}
\min_{u(t), t \in [\hatt_i,\hatt_{i+1}], u \in L^2} & \int_{\hatt_i + \invLa}^{\min\{\hatt_{i+1}+\invLa,T\}} OF_{\text{pen}}(Y_s,\hatt_i,y_{\hatt_i},y(s)) ds \notag \\
\ \text{subject to} \ \eqref{eq:densityOUP} \ \text{and} & \ z_t + \lambda z_x = 0, \quad z(0,t) = u(t), \quad z(x,\hatt_i) =  z_{\text{old}}(x,\hatt_i),\notag \\
& x \in (0,1), t \in [\hatt_i,\min	\{\hatt_{i+1}+\invLa,T\}], \label{eq:SOC-Prob_update}
\end{align}
where $z_{\text{old}}(x,\hatt_i)$ denotes the state of the production line at update time $\hatt_i$ ensuring that the SOC problems on the subintervals are correctly linked to each other.

Note that the usage of the demand density \eqref{eq:densityOUP} enables us to tackle both the SOC problem \eqref{eq:SOCmod} and the subproblems \eqref{eq:SOC-Prob_update} with methods from deterministic optimization, which will be done in the next section.

\section{A case study: The effect of an undersupply penalty}
\label{sec:numSim}
In this section, we numerically analyze the effect of an undersupply penalty for different intensities $\alpha$ for control methods $CM1$ and $CM2$. Using the reformulations \eqref{eq:SOCmod} and \eqref{eq:SOC-Prob_update} of the original SOC problem \eqref{eq:SOC}, the nonlinear optimization solver \textit{fmincon} from MATLAB R2015b\footnote{\url{https://de.mathworks.com/help/optim/ug/fmincon.html}, last checked: Sept 21, 2018} is applicable.

A left-sided Upwind scheme \cite{LeVeque.1990}, i.e. $\frac{z(x_j,\tau_{i+1})-z(x_j,\tau_i)}{\Delta \tau} + \lambda \frac{z(x_j,\tau_i)-z(x_{j-1},\tau_i)}{\Delta x}=0$, is chosen to discretize the linear advection equation. The applied step sizes $\Delta x=0.1$, and $\Delta \tau=\nicefrac{\Delta x}{\lambda}$ fulfill the CFL-condition.
For our numerical simulations, we use $10^3$ Monte Carlo repetitions with the following parameter setting for the demand process: $T=1$, $\lambda=4$, $\mu(t)=2 + 3 \cdot \sin(2\pi t)$, $\kappa=3$, $\sigma=2$, $y_0=1$.

\begin{figure}[h!]
	\subfloat[\ $\alpha=1$\label{fig:outputMMOUPwithUp_lambda4_initial1_2plus3sin2pit_kappa3_sigma2_nu5_gamma0_T1_alpha1_up5_MC1000path2}]{\includegraphics[width=0.49\textwidth]{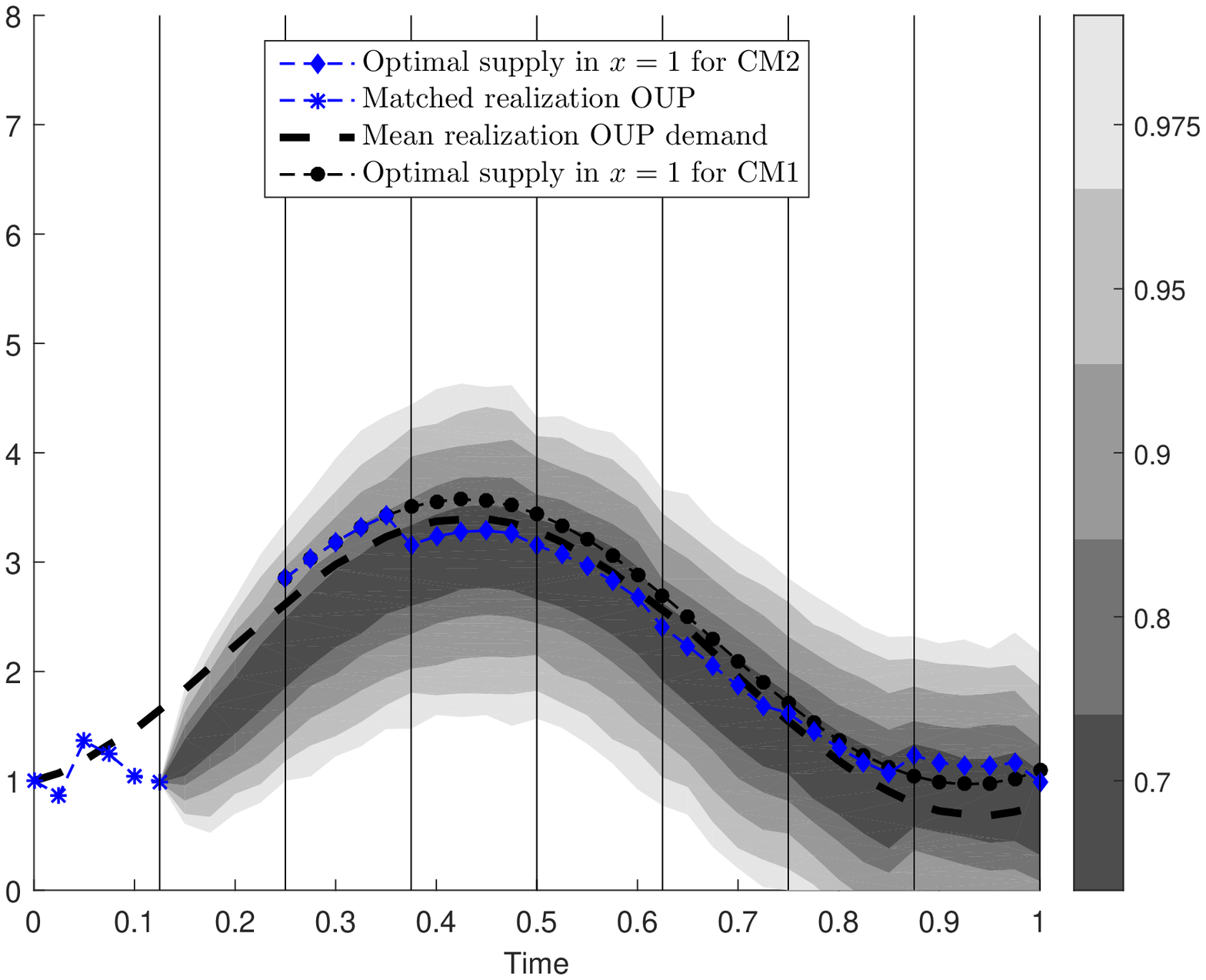}}\hfill
	\subfloat[\ $\alpha=3$\label{fig:outputMMOUPwithUp_lambda4_initial1_2plus3sin2pit_kappa3_sigma2_nu5_gamma0_T1_alpha3_up5_MC1000path2}]{\includegraphics[width=0.49\textwidth]{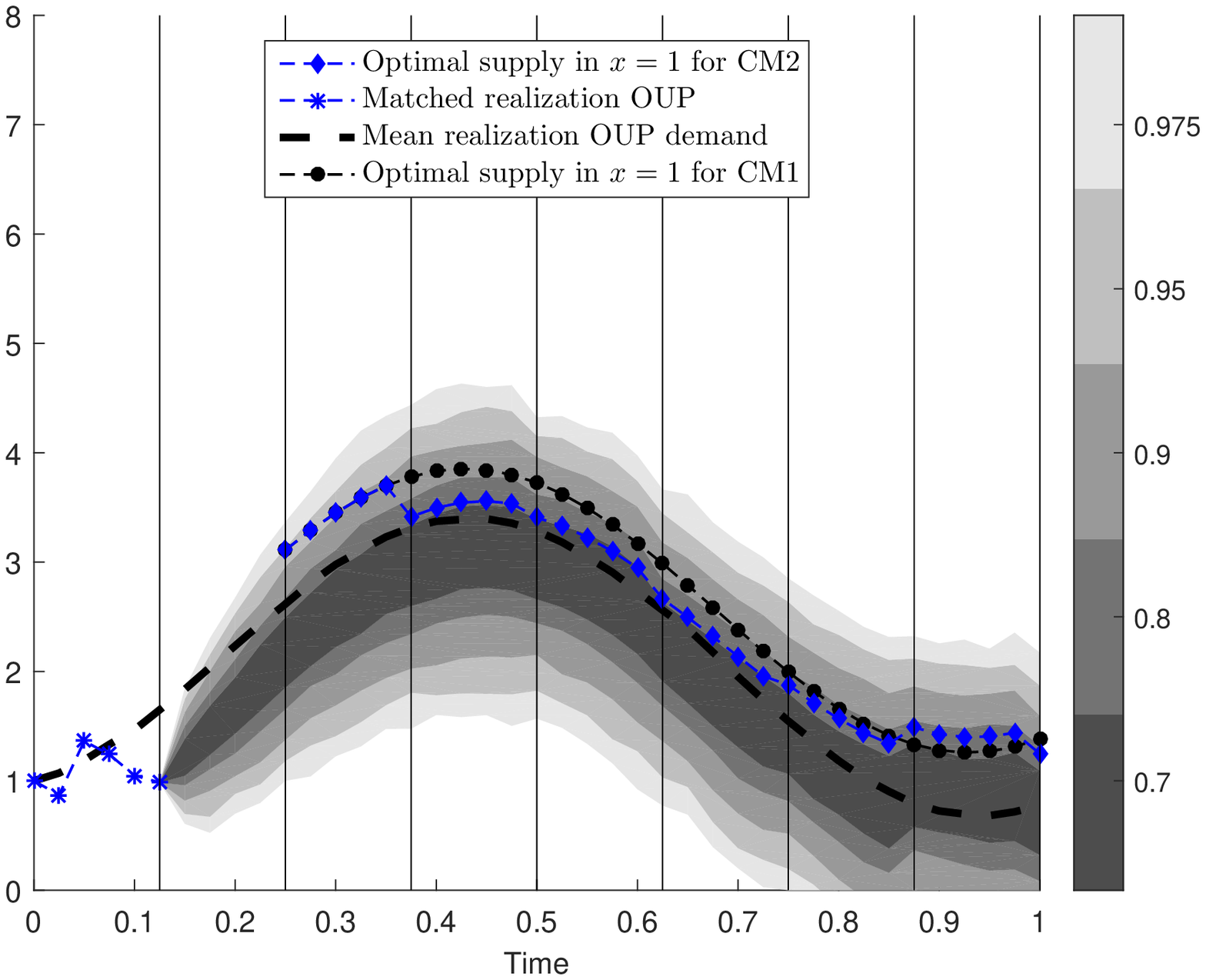}}
	\caption{Influence of penalty parameter $\alpha$ on optimal output}
	\label{fig:outputMMOUPwithUp_lambda4_initial1_2plus3sin2pit_kappa3_sigma2_nu5_gamma0_T1_up5_MC1000path2}
\end{figure}

In Figure \ref{fig:outputMMOUPwithUp_lambda4_initial1_2plus3sin2pit_kappa3_sigma2_nu5_gamma0_T1_up5_MC1000path2}, we are concerned with the influence of the penalty parameter $\alpha$ on the optimal output $y(t)$ for control methods CM1 and CM2. Thereby, we depict the updated confidence levels of the demand process in grey scale, the original mean realization of the demand (dashed line), the optimal CM1-output (dotted line), the optimal CM2-output (line marked by diamonds), and the tracked demand path until the first update time (line with asterisks). The vertical lines indicate the update times.
For both control methods, the penalty leads to an output above the (updated) mean demand. However, the CM1-output follows well the course of the original mean demand, and the CM2-output lies well within the upper part of the updated confidence intervals. Consistent with our intuition, a higher penalty parameter $\alpha$ leads to an output higher above the (updated) mean demand.
%Furthermore, we see that updates help to reduce oversupply. As the demand process is tracked to be lower than the mean demand level at the first update time (0.125), for example for $\alpha =1$, less output needs to be provided to still be in the upper part of the $70\%$-range of the demand process.

\begin{figure}[h!]
	\subfloat[\ Number of undersupply cases\label{fig:undersuppOUPwithUp_lambda4_initial1_2plus3sin2pit_kappa3_sigma2_nu5_gamma0_T1_up1_MC1000alpha13}]{\includegraphics[width=0.49\textwidth]{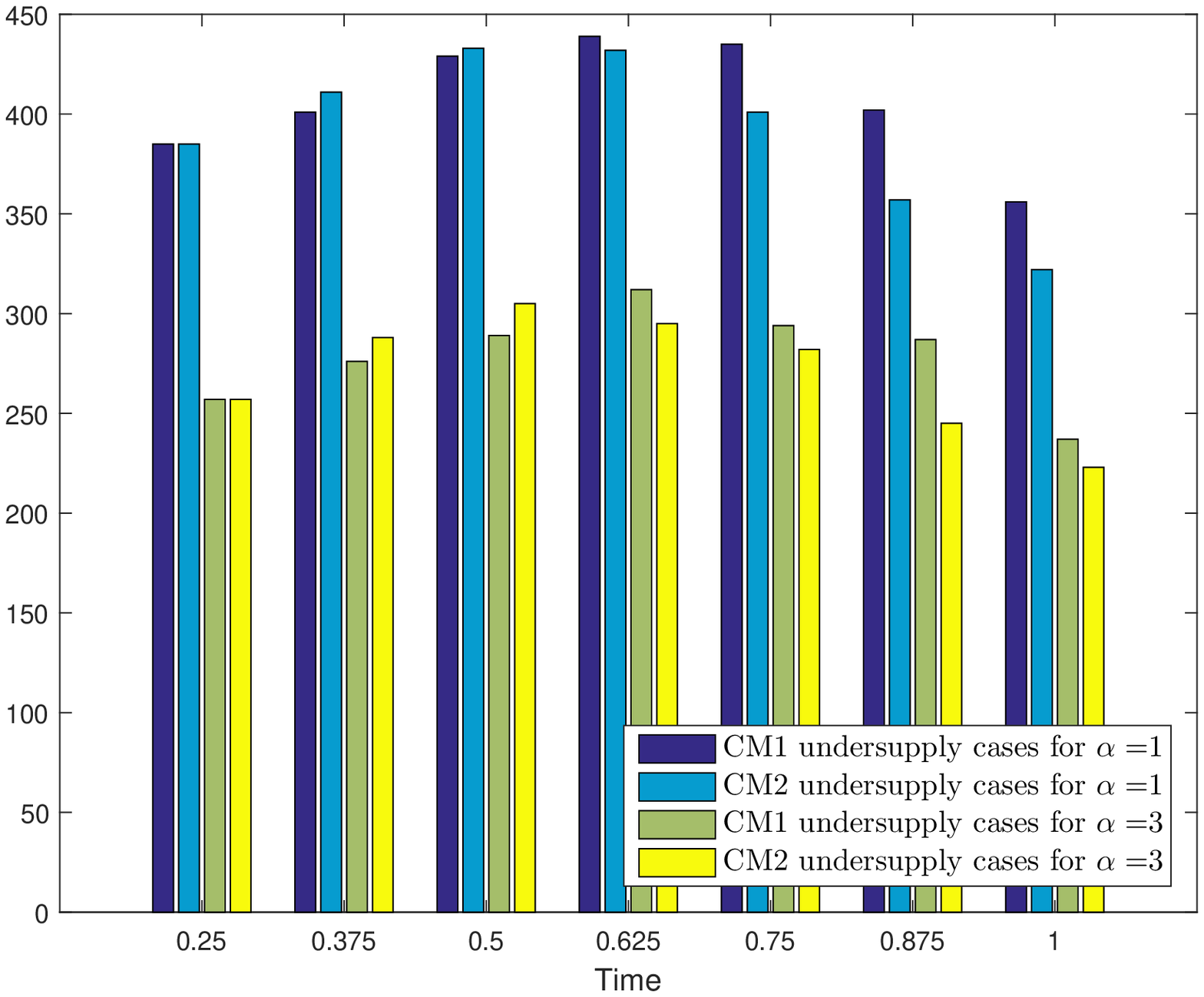}}\hfill
	\subfloat[\ Average undersupply \label{fig:notEnoughOUP_lambda4_initial1_2plus3sin2pit_kappa3_sigma2_nu5_gamma0_T1_up5_MC1000alpha13}]{\includegraphics[width=0.49\textwidth]{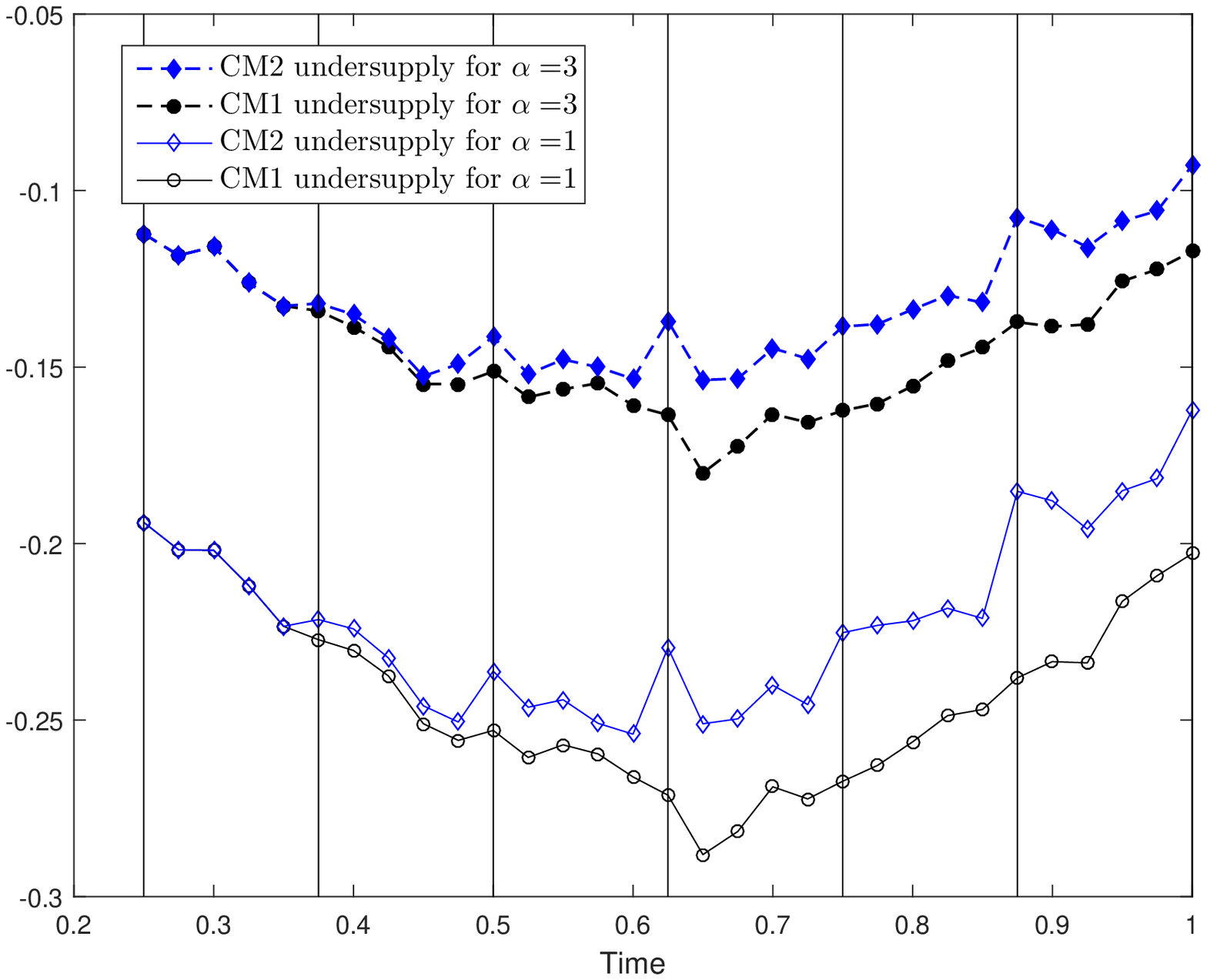}}
	\caption{Number of undersupply cases and average undersupply}
	\label{fig:underSupply}
\end{figure}

In a next step, we want to quantify the number of undersupply cases, i.e., for each point in time, we count how many of the $10^3$ simulated paths lie above the output (see Figure \ref{fig:undersuppOUPwithUp_lambda4_initial1_2plus3sin2pit_kappa3_sigma2_nu5_gamma0_T1_up1_MC1000alpha13}). By increasing the penalty parameter from $\alpha=1$ to $\alpha=3$, we are able to drastically reduce the number of undersupply cases.
Based on this information, it is not clear whether CM2 is preferable over CM1 or not. Note that deciding on an undersupply is a binary decision. However, in the objective function, the height of the deviation plays an important role. As there is a tradeoff between not realizing an undersupply but at the same time providing an adequate tracking of the demand, it might pay off to accept a small undersupply.
However, with respect to the average undersupply, Figure \ref{fig:notEnoughOUP_lambda4_initial1_2plus3sin2pit_kappa3_sigma2_nu5_gamma0_T1_up5_MC1000alpha13} shows that updates help to enhance the performance. To see this, at each point in time, we consider only those realizations where an undersupply occurs and plot the average height of the realized undersupply. The average undersupply for CM2 (lines marked by diamonds) is less or equal to the average undersupply for CM1 (dotted line). Furthermore, there is less average undersupply for a higher penalty parameter.
Finally, we can conclude that the introduction of a penalty parameter in the cost function leads to a reduction of both the undersupply cases as well as the average height of the undersupply.

\begin{acknowledgement}
The authors are grateful for the support of the German Research Foundation (DFG) within the project ``\textit{Novel models and control for networked pro\-blems: from discrete event to continuous dynamics}'' (GO1920/4-1) and the BMBF within the project ENets.
\end{acknowledgement}

\end{document}